%% file: ms.tex
\documentclass{birkjour_t2}

\usepackage{color}
\usepackage[dvipsnames]{xcolor}
\usepackage{listings}
\input{listings-python.tex}

\usepackage[hide]{ed} 
\usepackage{booktabs}
\usepackage{caption}
\usepackage{subfig}
\usepackage{wrapfig}

\usepackage{tikz}
\newcommand*\circled[1]{\tikz[baseline=(char.base)]{
            \node[shape=circle,draw,inner sep=1.3pt] (char) {#1};}}

\usepackage[hyperref,backend=bibtex,style=numeric]{biblatex}
\renewbibmacro*{event+venue+date}{}
\renewbibmacro*{doi+eprint+url}{%
  \iftoggle{bbx:doi}
    {\printfield{doi}\iffieldundef{doi}{}{\clearfield{url}}}
    {}%
  \newunit\newblock
  \iftoggle{bbx:eprint}
    {\usebibmacro{eprint}}
    {}%
  \newunit\newblock
  \iftoggle{bbx:url}
    {\usebibmacro{url+urldate}}
    {}}

\newcommand{\noop}[1]{}

\makeatletter
\newcommand\namedlabel[1]{#1\def\@currentlabel{#1}}

\def\blx@maxline{77}
\makeatother

\begin{document}

%
%
%
%
%
%
%
%
%

%
 \newtheorem{thm}{Theorem}[section]
 \theoremstyle{definition}
 \newtheorem{defn}[thm]{Definition}
\newtheorem{uc}{Use case}
\newtheorem{ex}{Example}
\newtheorem{prob}{Obstacle}

\title[DiscreteZOO]{DiscreteZOO: a Fingerprint Database of Discrete Objects}

\author[Ber\v{c}i\v{c}]{Katja Ber\v{c}i\v{c}}
\address{FAU Erlangen-N{\"u}rnberg}
\address{https://orcid.org/0000-0002-6678-8975}
\email{katja.bercic@fau.de}
\author[Vidali]{Jano\v{s} Vidali}
\address{Faculty of Mathematics and Physics,
University of Ljubljana, Slovenia}
\address{Institute of Mathematics, Physics and Mechanics,
Ljubljana, Slovenia}
\address{https://orcid.org/0000-0001-8061-9169}
\email{janos.vidali@fmf.uni-lj.si}


\date{November 30, 2018}


\begin{abstract}
In this paper, we present DiscreteZOO, a project 
which illustrates some of the possibilities for 
computer-supported management of collections
of finite combinatorial (discrete) objects, in particular
graphs with a high degree of symmetry.
DiscreteZOO encompasses a data repository, a website and 
a SageMath Package.
\end{abstract}

\maketitle


\section{Introduction}
\label{s:intro}

The fairly new discipline of Mathematical Knowledge Management 
\emph{``studies the possibility of computer-supporting and even automating 
the representation, cataloguing, retrieval, refactoring, plausibilisation, change propagation 
and even application of mathematical knowledge''}~\cite{Kohlhase:mkmtobbtg14}.
The discipline concerns itself with all forms of mathematical knowledge,
including definitions, theorems, proofs and data.
In mathematics, we commonly use the complementary strengths of 
humans and computers in creation of new knowledge,
while computer support for managing existing
mathematical knowledge is more rare.

\begin{wrapfigure}[15]{R}{3.5cm}
\footnotesize
\vspace*{-1.1em}
$G = (V, E)$, where $E = \{ \{i, j\}: i, j \in V \}$
\vspace{0.3cm}

\includegraphics[width=3cm,keepaspectratio]{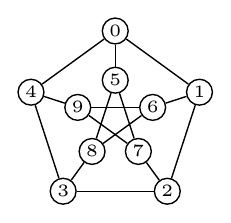}
\caption{An abstract object (top) and a concrete one (bottom).}
\label{figure:concrete-abstract}
\end{wrapfigure}

This is also the case when it comes to the management 
of a specific kind of mathematical data: collections of \textbf{concrete} 
(as opposed to abstract) mathematical objects
(Figure~\ref{figure:concrete-abstract}).
The DiscreteZOO project presented in this paper explores the possibilities 
for computer-supported management of collections of concrete examples
of finite combinatorial (discrete) objects.
We use collections of graphs with a high degree of symmetry as a
case study towards databases of a broader range of mathematical objects.

Mathematicians have a history of constructing catalogues, tables and 
similarly organised collections of mathematical objects.
In a thread started by Gordon Royle~\cite{RoyleMO}, 
the MathOverflow community collected several examples that predate computers.
Such collections can be used to better understand the objects in focus
(such as graphs with a high degree of symmetry) 
by studying small examples.
In particular, they can be computed with, 
used to search for patterns (or counter examples), 
to uncover connections, and to experiment with.
In combinatorics, enumeration of objects is sometimes done 
by actually constructing them and in some cases, 
it makes sense to store the objects as a catalogue.
These collections are also where mathematics can fruitfully intersect with other disciplines.

Broadly speaking, a fingerprint is something
that identifies something else 
(such as a person or data), often uniquely.
They are widely used in science and technology,
from forensics to genetics, biochemistry and computing%
~\cite{WikiFingerprint}.
In the context of mathematical knowledge, it is possible to 
organise fingerprints into a searchable database,
i.e., a \textbf{fingerprint database},
to support querying for mathematical results~\cite{BilleyTenner}.
Such databases make the data more accessible through querying and searching,
and thus enable mathematicians to be more effective at their work.
Last but not the least,
it is easier for a collection of mathematical objects
to get exposure when it is organised into a searchable database.

\subsection{Overview}

In Section~\ref{s:soa},
we sketch the state of art where it comes 
to collections of concrete mathematical objects.
In Section~\ref{s:req},
we present some use cases and combine them with the state of the art
to identify certain requirements we consider necessary
to bridge the experience gap  between the tools that are currently offered
and what the users (particularly mathematicians) need.
In Section~\ref{s:desc}, we describe the DiscreteZOO project,
paying particular attention to the structure of information and 
design decisions.%
We conclude in Section~\ref{s:fw} with an outline of future work.

\bigskip
\paragraph*{Acknowledgements}%
The work by the first author reported here is supported in part
by the OpenDreamKit Horizon 2020 European Research Infrastructures project
(\#676541).
The first author did the main portion of her work on the project
while a postdoctoral scholar at the Center for Mathematical Sciences
of the National Autonomous University of Mexico, Campus Morelia.
The work by the second author reported here
is supported by the Slovenian Research Agency
(research program P1-0285 and project J1-8130).
The first author would also like to thank Michael Kohlhase
for his encouragement, his many insightful comments,
and for pointing out the need to include a checksum in a CID.

\section{State of the Art}
\label{s:soa}

We use the expression \textbf{collection of mathematical objects} (CMO)
to describe aggregations of concrete mathematical objects.
In other contexts, these have been called
atlases, catalogues, censuses, databases, directories, encyclopedias, and lists.
Some of these are associated with some structure or features.
Encyclopedias might be expected to have one interface page per object.
A census is usually a collection that is in some way complete.
For CMOs, we want to presuppose neither an interface
nor any kind of underlying system.
The only requirement we impose is a collection
only contains concrete objects of one type.
(such as collections of integer sequences or graphs).
We do not require any degree of digitalisation for collections
(and so, the \emph{Atlas of Graphs}~\cite{AtlasOfGraphs} contains a collection of graphs)
nor do we want to bind collections of mathematical objects to any kind of infrastructure.
We will use the word \textbf{database} whenever the collection
has infrastructure that supports querying.

We know~\cite{CMO} of at least nine
searchable databases of mathematical objects,
which mostly address the needs of specific communities:
\emph{ATLAS of Finite Group Representations (Version 3)}~\cite{GroupRepresentations},
\emph{The Database of Permutation Pattern Avoidance}~\cite{PermutationPatterns}
\emph{Encyclopedia of Graphs}~\cite{GReGAS},
\emph{FindStat}~\cite{findstatlocal}, 
the \emph{House of Graphs}~\cite{HOG},
the \emph{LMFDB}~\cite{lmfdb:on},
the OEIS~\cite{oeis},
the \emph{QaoS}~\cite{QaoS},
the \emph{Small Graph Database}~\cite{SmallGraph}
and the \emph{SymbolicData Project}~\cite{SymbolicData}.
We found many more CMOs with fewer features,
for example five collections of abstract polytopes,
and, depending how one counts,
several collections of lattices (at least three), graphs (at least ten) 
and other, mostly combinatorial, objects (at least ten).
This is not at all surprising. 
On one hand, the accessibility and advancement of technology
made it possible for more and more researchers to produce
collections of mathematical objects.
On the other hand, development of nontrivial infrastructure
necessary for a searchable database goes beyond 
a typical mathematician's resources, skillset or interests.

There are some theoretical obstacles one encounters
when making mathematical knowledge more accessible
that are also true when working with CMOs.

\begin{prob}
Most mathematical results or objects do not have 
a natural order or structure that could be used 
as an intuitive basis for browsing.
\end{prob}

The standard way to approach this is to support querying.
Billey and Tenner proposed~\cite{BilleyTenner} \textbf{fingerprints},
or canonical representations of mathematical results,
as a means to querying databases of theorems -- 
an idea that also makes sense in the context of CMOs.
Ideally, all fingerprints would be small, language-independent
and would uniquely determine the fingerprintees.
Finding suitable fingerprints is a domain-specific mathematical problem
that does not always produce a result.

\begin{prob}
Many mathematical results or objects do not have 
a representation that is small, language-independent and canonical
(uniquely determining).
\end{prob}

Notable examples of objects that do not have a canonical form are graphs.
They do, however, have a form that comes close.
A \textbf{canonical form} or labelling~\cite{Canonical} of a graph $G$ is a labelled graph
$\operatorname{Canon}(G) \simeq G$ such that
that for every other graph $H$, $H \simeq G$ if and only if
$\operatorname{Canon}(H) = \operatorname{Canon}(G)$,
when both $\operatorname{Canon}(H)$ and $\operatorname{Canon}(G)$
can be computed.
The labelling obtained depends on the algorithm that was used to compute it.
Nauty~\cite{NautyTraces} and Bliss~\cite{Bliss} are two libraries 
commonly used to produce canonical forms of graphs.

The existence of a natural fingerprint likely played a part in the success of 
the \emph{On-Line Encyclopedia of Integer Sequences} (OEIS)~\cite{oeis},
perhaps the most famous example of a database of mathematical results.
The OEIS is a searchable and collaborative database of integer sequences,
which serve as fingerprints for their associated entries.  
The indexing is based on small, language-independent, and canonical data:
the first elements of a sequence.

Even if most CMOs are only weakly supported by computers,
they are already being used to 
manually find patterns (by writing programs) in the CMOs
and look up properties of mathematical objects.
A searchable database would make working with them more efficient and would
open up additional options to researchers.
Indeed, the initial motivation behind the DiscreteZOO project was
to provide a home for the partially overlapping censuses of 
graphs with a high degree of symmetry~\cite{FosterCurrent,VTCensus,CVTCensus},
even if it quickly became clear that there was no reason
to restrict it to being a database for symmetric graphs.

It is important to be able to refer to a mathematical object 
by a short unique identifier, for example in a paper,
and some of the searchable databases we just mentioned
already provide citable indexes.
The sequences in the OEIS are citable via 
their unique identifier (such as \texttt{A000055}) 
Databases like
\emph{The Database of Permutation Pattern Avoidance}~%
\cite{PermutationPatterns}
and \emph{FindStat}~\cite{findstatlocal}
use similar citable unique identifiers.
A a system of Math Object identifiers (MOI)
has been proposed in~\cite{Kohlhase:moitrdm17},
which tries to solve the CID problem
for all kinds of mathematical knowledge.

While some of the more feature-rich databases of 
mathematical objects accept new contributions,
there is currently no other simple way to organise a 
collection of mathematical objects into a searchable database.
On the other hand, there are some helpful features that 
are not implemented even by the existing searchable databases.
The LMFDB is a good example of a database with an interface
that goes as far as it can, but does not have an API
for programmatic access.

As mathematicians find themselves working with data more frequently
and as the industry standards for interacting with data move forward,
the current state of the art in becomes insufficient.

\section{Requirements Analysis}
\label{s:req}

There is a need in the mathematical community for 
a platform or framework that would 
provide infrastructural support for CMOs.
We believe this is possible, because while every CMO
necessarily has its own peculiarities,
they have enough in common for a reusable infrastructure to make sense.
This infrastructure could then provide common features to all databases
while being flexible enough to accommodate any database-specific features.
Such a platform would not need to implement
many or complex features, as most CMOs would benefit
even from just a very basic search functionality.

In the rest of this section, we will first describe some use cases
we use as guidelines in the design of the DiscreteZOO system.
In the second part of the section, we will describe 
some requirements for an ideal database of mathematical objects.
Most of the features derived from these requirements could be 
provided by a general purpose platform or framework for CMOs.

\begin{uc}[Classification of cubic vertex-transitive partial cubes]
\label{ex:partcubes}
Tilen is interested in cubic partial cubes that are also vertex-transitive%
\footnote{Tilen only used the original census for his research~\cite{Marc}.}.
The census of cubic vertex-transitive graphs
on up to $1280$ vertices~\cite{CVTCensus}
contains (among others) all the graphs he is interested in
on up to $1280$ vertices.
First, he needs to filter out those which are partial cubes.
He would like to know whether
there are any vertex-transitive cubic partial cubes
in addition to the infinite family of even prisms ($K_2 \square C_{2n}$).
Finally, he would like to compute isometric embeddings into hypercubes.

The census is published as Magma code~\cite{Magma},
and the only way to do what he wants (apart from DiscreteZOO)
is to program everything in a computer algebra system.
If he does not have access to Magma,
he needs to process the code
so that it is readable by another computer algebra system.
We will show how the task looks like with DiscreteZOO
in Sections~\ref{ss:website} and~\ref{ss:sage}.
\end{uc}

\begin{uc}[Classification of cubic vertex-transitive graphs with small girth]
\label{ex:smallgirth}
Primo\v{z} wants to work on a classification
of cubic vertex-transitive graphs with small girth~\cite{PotocnikVidali}.
To get an idea on how to proceed,
he first wants to know how many graphs from 
the census~\cite{CVTCensus} there are for each girth.
Then he wants to check some additional properties 
for each graph.
\end{uc}

\begin{uc}[Abstract polytope lookup] \label{ex:polytope}
Alice has recently come across an interesting abstract polytope.
She knows it is self-dual and self-Petrial with the Schl{\"a}fli symbol $\{12, 12\}$,
but does not remember its name.
She would like to look it up quickly in a database by filtering 
using the properties she remembers.
\end{uc}

\begin{uc}[Automorphism groups downloads] \label{ex:groups}
Bob wants to test a subgroup condition for automorphism groups of certain maps.
He would like to use a database of maps (or equivalently, maniplexes of rank 3)
to obtain a list of automorphism groups
which he could then test with a GAP function he had already written.
\end{uc}

Use cases~\ref{ex:partcubes} and~\ref{ex:smallgirth}
are already fully supported in DiscreteZOO.
Use case~\ref{ex:polytope} would require 
supporting properties which are numeric arrays 
(for the Schl{\"a}fli symbol property).
Use case~\ref{ex:groups} often appears in 
the area of symmetries of discrete objects,
where researchers are often
interested in properties of automorphism groups
of objects which exhibit some specific property. 
It is beyond the scope of DiscreteZOO
to provide such an in-depth analysis of automorphism groups,
especially since computer algebra systems serve that purpose well.
However, we believe that generating CAS code for 
lists of automorphism groups would be a valuable time-saving feature.

\subsection{Requirements}

Billey and Tenner~\cite{BilleyTenner} identified several properties 
important for fingerprint databases of theorems.
Namely, such a resource should be a
\emph{``searchable, collaborative database of citable mathematical results, 
indexed by small, language-independent, and canonical data''}.
Given the rate of growth of mathematical knowledge,
we believe these features to be (or become) critical for any mathematical knowledge resource,
like databases of definitions, objects, theorems and proofs.
We added some additional requirements we believe to be important
in the following list.

\begin{description}
\item[\namedlabel{R1}: Searchable\label{re:searchable}]
The database should have at least
a basic search (filter) functionality for objects
\item[\namedlabel{R2}: Collaborative\label{re:collaborative}]
The database should be collaborative.
It needs to be easy for people to contribute and 
it must be possible for any data to find out who contributed it and when.
In particular,
\begin{enumerate}
\item \label{re:declarative} adding a new collection should ideally be a declarative task
as opposed to a programming one --
it should suffice to describe the structure of the collection,
\item changes such as adding new kinds of objects, new collections, new properties, 
or updates to existing values need to be tracked by the database.
\end{enumerate}
\item[\namedlabel{R3}: Citable\label{re:citable}] The records in the database need to be citable.
\item[R4: User-friendly] The database should provide suitable interfaces.  This means
\begin{enumerate}
\item simple, intuitive and easily accessible interfaces for casual users
(it should avoid exposing the casual user to low level interfaces,
such as Python, SageMath, SQL, \ldots),
\item making usage easy or easier for power users.
\end{enumerate}
\item[\namedlabel{R5}: Self-explaining\label{re:self-explaining}]%
The database interfaces should make accessing
definitions of concepts and further information easily available.
\item[\namedlabel{R6}: Interoperable\label{re:interoperable}]%
Other systems (like databases and computer algebra systems) 
should be able to interact with the database.
This involves having well defined APIs.
\item[\namedlabel{R7}: Coverage information\label{re:coverage-information}]%
The database interfaces should make explicit
the assumptions made by the system and collections as well as providing 
information on the completeness of the search results.
For example for graphs:
\emph{``the search results are complete for your 
search parameters up to order 511''}.
\item[R8: Non-redundant] The objects should be stored up to isomorphism.
\item[\namedlabel{R9}: Decentralised\label{re:decentralised}]%
The database should not rely on the existence of a central authority.
An entry in a local copy of the database
should be easily distributable to other copies
without any third-party intervention.
\end{description}

Another reason for a decentralised database is the
following scenario.
A researcher may encounter an object
that is not yet in the database;
another researcher may then reproduce her work (or use his own methods)
and easily verify that they have obtained the same object
by comparing the identifier.
Once the object is in the database,
the identifier may then be used to quickly access the object
and its properties.

\section{DiscreteZOO, a First Step Towards a Generic CMO Framework}
\label{s:desc}

DiscreteZOO is an attempt to provide a generic platform with 
few features that are as general as possible,
taking as a use case collections of symmetric graphs and related objects.
We chose to focus on features that are simplest to implement 
and that have the greatest usability.
We tried to keep the infrastructure as ``object type agnostic''.
That is, we avoided features specific to the object type
(except for actual database encodings).
An object types comes with a set of properties 
(mathematical invariants) that apply to objects of this type
which are then used for searching and filtering.
Here, we focused on property types that are easiest (and, in fact, easy) 
to implement: Boolean and numeric properties.

There are two main user interfaces: the website and the SageMath package.
Ideas involved in the design process~\cite{DesignSystems,EverydayThings}
were helpful in the preparation of the interfaces, 
especially so for the website.
The interfaces implement a simple query language which we believe
to be a flexible and general solution for accessing results in CMOs.

New data come into DiscreteZOO 
through the data repository, 
while the class specification repository 
describes the objects and their types.
The need for these two repositories
arose from the fact that this information
is shared between different applications,
with different needs, and
written in multiple programming languages.
Therefore, the data has to be represented
in a programming language-independent manner.
We chose the JSON format for its widespread support,
and for being relatively human readable.
Additionally, the usual line-based formatting
means that such a format is well-suited for 
tracking changes with Git, 
which also addresses~\ref{re:collaborative}.

The layout of the data repository also closely follows the one used by Git,
which stores its objects in files whose paths determine their hashes.
Unlike Git,
which essentially uses a content-addressable filesystem~\cite{Git},
the DiscreteZOO data repository is not entirely content-addressable.
Instead, it could be described as {\em object-addressable} --
for each object, the path to the file describing it (or a symbolic link to it)
can be computed from the object itself.
However, this file may still change (and keep the same path)
if, say, new properties for this object are added.

\subsection{Information Structure}

DiscreteZOO aims to provide a platform for CMOs with a structure that is particularly simple:
ones than can be represented as tables, where each row represents information about one object
and where each column represents a mathematical property (invariant) of objects.
The actual database schema might be more complex for implementation purposes.
Moreover, we require for now that these properties have
either Boolean values (a graph is either Hamiltonian or not)
or numeric (integer or real) values,
small enough to fit standard database types.

In general, it would be possible to store any objects that have
a reasonably small encoding as text or a binary blob.
Each of the object types comes with a set of properties
that are relevant for it,
e.g., the properties that make sense for symmetric graphs
are not the same as the ones that make sense for maniplexes.

A typical DiscreteZOO database record
describes a single discrete object, 
several of its mathematical properties,
indexes in collections,
and any number of human-readable aliases and descriptions.
Most of the properties are mathematical invariants.
Every object also has a \textbf{GUID} (globally unique identifier)
and a \textbf{CID} (citable identifier).
For more detail on partial fingerprints, GUIDs and CIDs,
see~\ref{sss:fingerprints}.
All objects included so far also have 
a shareable encoding consisting of printable characters.  
Such an encoding makes it possible to transfer an object between software tools.

For example, the record representing the Petersen graph
(Table~\ref{tb:petersen})
is marked as a graph, vertex-transitive graph
and a cubic vertex-transitive graph
and has all the precomputed properties relevant to these classes.

\begin{table}
\smaller
\begin{tabular}{lcl}

\multicolumn{3}{l}{
\begin{tabular}{ll}
\textbf{GUID} & \texttt{c74c6028a25a65a6189db885bdaa11aa1c1f2a861f4c0d0e8efd6a4ec786e0a5} \\
\textbf{CID} & \texttt{Zc74c-6028-a25a+8} \\
\textbf{data} & \texttt{:IeIKqPD?hgAH?G\textasciitilde} \\
\textbf{aliases} & Petersen graph \\
\end{tabular}
\vspace{0.5cm}
} \\

\begin{tabular}[t]{ll}
\multicolumn{2}{l}{Class \emph{graph}} \\
\midrule
is arc-transitive & true \\
is bipartite & false \\
is Cayley & false \\
is distance-regular & true \\
is distance-transitive & true \\
is edge-transitive & true \\
is Eulerian & false \\
is Hamiltonian & false \\
is overfull & false \\
is partial cube & false \\
is split & false \\
is strongly regular & true \\
\midrule
clique number & $2$ \\
connected components number & $1$ \\
diameter & $2$ \\
girth & $5$ \\
odd girth & $5$ \\
order & $10$ \\
size & $15$ \\
triangles count & $0$ \\
\bottomrule
\end{tabular}

&
\hspace{0.5cm}
&

\begin{tabular}[t]{ll}
\multicolumn{2}{l}{Class \emph{vertex-transitive graph}} \\
\midrule
Index in~\cite{VTCensus} & $(10, 5)$ \\
\bottomrule
\multicolumn{2}{l}{\vspace{1cm}} \\
\multicolumn{2}{l}{Class \emph{cubic vertex-transitive graph}} \\
\midrule
Index in~\cite{FosterCurrent} & $(10, 1)$ \\
Index in~\cite{CVTCensus} & $(10, 3)$ \\
\midrule
is M\"obius ladder & false \\
is prism & false \\
is SPX & false \\ 
\bottomrule
\end{tabular}
\\
\end{tabular}
\vspace{0.5cm}
\caption{The Petersen graph}
\label{tb:petersen}
\end{table}

\subsubsection{Fingerprints}
\label{sss:fingerprints}

Billey and Tenner~\cite{BilleyTenner} observed the natural conflict between 
keeping fingerprints small on one hand and their canonicity on the other.
If we want to be able to use fingerprints for querying,
they should also be readable or writable by humans
(for example, for entry in a search input field on a website).

Mathematical invariants, such as valency or girth for graphs,
behave as fingerprints in some ways
and we believe it makes sense to use them as such.
They are examples of data that are 
small, language independent, 
but do not necessarily uniquely identify an object:
for example, there are infinitely many arc-transitive graphs of valency $3$.
We propose the term \textbf{partial fingerprints}
to describe data that are small and language independent,
but do not uniquely identify objects.
In DiscreteZOO, we use the Boolean and numeric properties
as partial fingerprints.

Ideally, the objects in DiscreteZOO have a GUID based on 
some canonical string representation.
If such a representation exists, we use the
$64$-character hexa\-decimal representation of its SHA-256 hash for the GUIDs.
Such GUIDs are fingerprints in the sense that they display a level of canonicity 
(two GUIDs obtained by the preferred algorithm are the same if and only if the objects are the same and 
it is possible to decode the complete object from them) 
and are reasonably small and language-independent,
but are not human friendly.
We plan to use randomly assigned GUIDs if we ever encounter 
an object type that has nothing approaching a canonical form.
The use of a cryptographic hash function
is a strong guarantee~\cite[\S9]{CryptoHandbook}
that no two objects are ever expected to be found
to have the same identifier.
The reason for choosing such a hash-based identifier
instead of a sequential one (as with other databases) is decentralisation (\ref{re:decentralised}).

For graphs, the canonical labellings are 
as close as it is currently possible to get to a canonical form.
We base the GUIDs for graph on the \texttt{graph6} and \texttt{sparse6}
graph formats~\cite{GraphFormats}.
It is possible to use the graph canonical labellings to
obtain canonical representations of other graph-like objects:
for maniplexes, we can just add the edge labelling.

\subsubsection{Citability}

Ideally, a database of mathematical objects provides both an easy way to cite its contents as well as 
an efficient way to compare its contents with another collection of mathematical objects.
These two objectives are, in some way, at cross purposes.
The former requires the citable identifier to be as short as possible, 
while the other most likely means we need to use some kind of hashes (large, to avoid conflicts).

The $64$-character strings of the GUIDs in DiscreteZOO
are infeasibly long for reproduction in a text intended to be read by humans.
We chose a standard abbreviation technique.
In the Git versioning system~\cite{Git},
the objects are identified by $40$-character hashes, 
but are usually referred to by simply taking the first $7$ characters.
In DiscreteZOO, the citable identifier is obtained
by taking the first $12$ hexadecimal digits of the hash
and using further characters when necessary to avoid conflicts.
A checksum hexadecimal digit,
obtained by taking the XOR of the digits included in the citable identifier,
is then added to the end.
For readability, the characters are split into groups of $4$ characters,
the letter \texttt{Z} is prepended,
and the checksum digit is separated by a plus sign:

\[
\texttt{123456789abcdef...} \quad \rightarrow \quad \texttt{Z1234-5678-9abc+c}.
\]

\bigskip
The DiscreteZOO citable identifier (\ref{re:citable})
is described in more detail in the project documentation~%
\cite[{\tt Documentation} repository]{ZooGroup}.

\subsubsection{Query Language}
\label{sss:query}

The query language supports conditions 
on Boolean and numeric property types, as described below.

\begin{defn}
Let $o$ be an object and let $B$ be a Boolean property of this object.
A \textbf{Boolean condition} is one of the following predicates.
\begin{itemize}
\item $o$ has property $B$,
\item $o$ does not have property $B$.
\end{itemize}
\end{defn}

\begin{defn}
\label{def:numeric}
Let $o$ be an object and let $N$ be a numeric property of this object
and let $n$ be a number and $R \in \{ =, \neq, <, \leq, >, \geq \}$.
A \textbf{numeric condition} is a predicate $N \ R \ n$.
\end{defn}

Examples of Boolean conditions for graphs
are ``is bipartite'', ``is not Hamiltonian'', etc.
Examples of numeric conditions for graphs
would be ``girth equal to $3$'', ``order at most $100$'', etc.
Currently, expressing a condition ``valency is prime'' is not possible
on the level of data 
(it is of course possible to check for it in the SageMath package).

It is only possible to query objects of the same type.
Both the website and the SageMath package support simple queries (\ref{re:searchable}),
consisting of conjunctions of the predicates described above.

\begin{figure}[h!]
Find all $o$ of type $T$ such that $P(o)$ holds for all $P \in \mathcal{C}$,
\\ where $\mathcal{C}$ is a set of Boolean and numeric conditions.
\end{figure}

The SageMath package allows arbitrary logical propositions
(e.g., conjunctions, disjunctions, exclusive or, etc.)
with the above predicates acting as atoms.
Arithmetic can also be performed on numeric properties
to derive new numeric conditions.
The results can also be ordered
according to some atomic or derived properties.

Let $\mathcal{P}$ be a logical proposition
containing only Boolean and numeric conditions
derived from the set of all properties of $T$.
Let us additionally define
$\mathcal{G}_1, \mathcal{G}_2, \dots, \mathcal{G}_k$
as properties derived only from Boolean and numeric conditions
from the set of all properties of $T$
(through logical manipulation of Boolean properties
and arithmetic manipulation of numeric properties).
The SageMath supports the following queries.

\begin{figure}[h!]
Find all $o$ of type $T$ such that $\mathcal{P}(o)$ holds,
in the lexicographic order of
$(\mathcal{G}_1(o), \mathcal{G}_2(o), \dots, \mathcal{G}_k(o))$.
\end{figure}

Additionally, besides iterating through the matches,
it also allows counting queries
(i.e., how many objects in the database satisfy the proposition).
These can also be further enriched by grouping,
so a single query can count the number of objects
sharing the same value of a property (either atomic or derived)
for all values appearing in the matches.

\begin{figure}[h!]
For each value
$(\mathcal{G}_1(o), \mathcal{G}_2(o), \dots, \mathcal{G}_k(o))$
of an object $o$ such that $\mathcal{P}(o)$ holds, \\
find the number of such objects.
\end{figure}

\begin{ex} \label{ex:cvtg5}
Find all graphs that are vertex-transitive,
of valency $3$ and girth at least $5$.
\end{ex}

\begin{ex} \label{ex:cvtgcount}
For each girth, count all cubic vertex-transitive graphs of diameter $5$.
\end{ex}

\subsection{Contents}

DiscreteZOO currently contains $212269$
graphs from the following three collections.
\begin{itemize}
\item The catalogue of vertex-transitive graphs on up to $31$ vertices
by Brendan McKay and Gordon Royle~\cite{VTCensus} ($100661$ graphs),
\item the census of cubic vertex-transitive graphs on up to $1280$ vertices
by Primo\v{z} Poto\v{c}nik, Pablo Spiga and Gabriel Verret~\cite{CVTCensus}
($111360$ graphs), and
\item the census of cubic arc-transitive graphs on up to $2048$ vertices
by Marston Conder~\cite{FosterCurrent} ($796$ graphs).
\end{itemize}

A collection of $75217$ maniplexes
is a work in progress by the first author.

Table~\ref{tb:props} shows the properties which are currently computed for
nearly all records and which do not have the same value for all of them.
There are roughly $60$ graph related properties in total, most of which
are Boolean, integer or real-valued.
We are also working on relations between objects
(graph-valued graph properties) such as truncation.

\begin{table}
\smaller
\centering
\subfloat
{%
\begin{tabular}{ll}
\toprule
Graph property &  \\
\midrule
is arc-transitive & Boolean \\
is bipartite & Boolean \\
is Cayley & Boolean \\
is distance-regular & Boolean \\
is distance-transitive & Boolean \\
is edge-transitive & Boolean \\
is Eulerian & Boolean \\
is Hamiltonian & Boolean \\
is M\"obius ladder & Boolean \\
is overfull & Boolean \\
is partial cube & Boolean \\
is prism & Boolean \\
is split & Boolean \\
is SPX & Boolean \\
is strongly regular & Boolean \\
\midrule
chromatic index & numeric \\
clique number & numeric \\
connected components number & numeric \\
diameter & numeric \\
girth & numeric \\
odd girth & numeric \\
order & numeric \\
size & numeric \\
triangles count & numeric \\
\bottomrule
\end{tabular}
}%
\qquad\qquad
\subfloat
{%
\begin{tabular}{ll}
\toprule
Maniplex property &  \\
\midrule
is polytope & Boolean \\
is regular & Boolean \\
\midrule
group order & numeric \\
number of orbits & numeric \\
rank & numeric \\
\bottomrule
\end{tabular}
}%
\caption{Properties of graphs and maniplexes in DiscreteZOO}
\label{tb:props}
\end{table}

\subsection{Class specifications and the data repository}
\label{ss:repos}

\begin{wrapfigure}{R}{0.45\textwidth}
\centering
\includegraphics[width=0.4\textwidth,keepaspectratio]{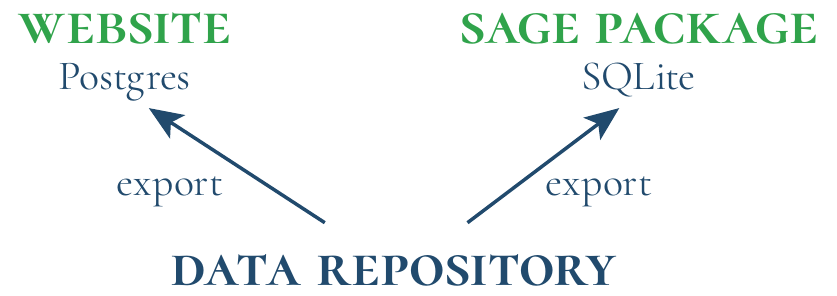}
\caption{Data repository and databases}
\end{wrapfigure}

The \textbf{class specifications}~%
\cite[{\tt Specifications} repository]{ZooGroup}
serve as a declarative (and to an extent, language-agnostic) way
(\ref{re:collaborative}.\ref{re:declarative}) to describe CMOs.
Together with the \textbf{data repository}~%
\cite[{\tt Data} repository]{ZooGroup},
they can be used to produce the databases
used by the website and the SageMath interface.
Each of them is available as a public repository on GitHub
containing JSON files with the relevant information.

The class specifications determine the possible classes of objects
that can be stored in the DiscreteZOO data repository.
The classes are organised in a hierarchical manner:
\pythoninline{ZooEntity} represents the top-level class,
and all other classes are descended from it.
Directly descending from \pythoninline{ZooEntity}
is the class \pythoninline{ZooObject},
which serves as a superclass for all mathematical objects in DiscreteZOO.
Classes not descending from \pythoninline{ZooObject}
are used to represent metadata --
for instance, there is a \pythoninline{Change} class
which is used to track changes to the database.

\begin{wrapfigure}{R}{0.45\textwidth}
\centering
\vspace*{-0.5em}
\includegraphics[width=0.4\textwidth,keepaspectratio]{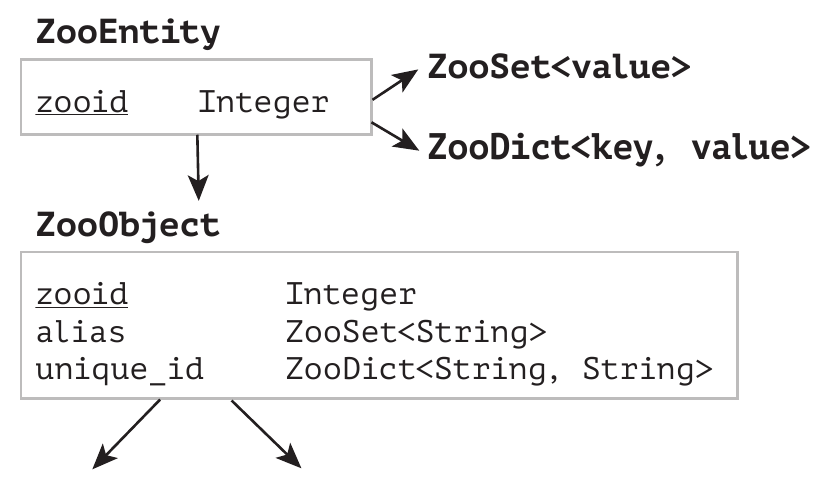}
\caption{Class hierarchy}
\end{wrapfigure}

Each class specification is a JSON file
storing an object which includes the following fields:
\begin{itemize}
\item \pythoninline{fields}:
an object mapping property names to their types
(one of \pythoninline{Integer}, \pythoninline{Rational},
\pythoninline{RealNumber}, \pythoninline{bool}, \pythoninline{str},
or a DiscreteZOO class);
\item \pythoninline{primary_key}:
the name of the property determining an instance
(i.e., it acts as the primary key in the database);
\item \pythoninline{condition}:
an object specifying the required values of properties;
\item \pythoninline{default}:
an object specifying default values of properties
(i.e., those which are determined by the conditions).
\end{itemize}

Note that the \pythoninline{fields} object
only contains the names of the properties which are specific to the class --
its instances will then additionally have
all properties specified by the superclasses.
The parent class is determined by the type of the property
which is specified as \pythoninline{primary_key}
(for \pythoninline{ZooEntity}, this is \pythoninline{Integer},
which is a builtin type).

For example, the \pythoninline{ZooObject} class has three properties:
\pythoninline{alias}, \pythoninline{unique_id}, and \pythoninline{zooid}.
The latter is also specified as the primary key,
and its type \pythoninline{ZooEntity} means that it acts as a reference
to the primary key of the \pythoninline{ZooEntity} class,
i.e., its only property \pythoninline{zooid}.
The first two properties have types deriving from the classes
\pythoninline{ZooSet} and \pythoninline{ZooDict}, respectively.
These two classes can be used to construct set and dictionary types,
and do not have a corresponding JSON specification.
The two properties thus provide
every mathematical object stored in DiscreteZOO
with a set of aliases and a dictionary of GUIDs --
hashes of a canonical representations of an object
as given by a chosen algorithm.

The data repository stores information about the objects in DiscreteZOO.
Each object is stored in a JSON file whose path is determined by its GUID:
the first level folder gives the algorithm used,
the second level folder is named for the first two nibbles
in the hexadecimal representation of the GUID,
and the filename corresponds to the remainder of the GUID.
Since the GUID is only unique up to the algorithm used to compute it,
an object may be represented by more than one file.
In this case, only one of them is an ordinary file,
while the rest are symbolic links pointing to the actual file.

\begin{figure}
\centering
\includegraphics[width=0.8\textwidth,keepaspectratio]{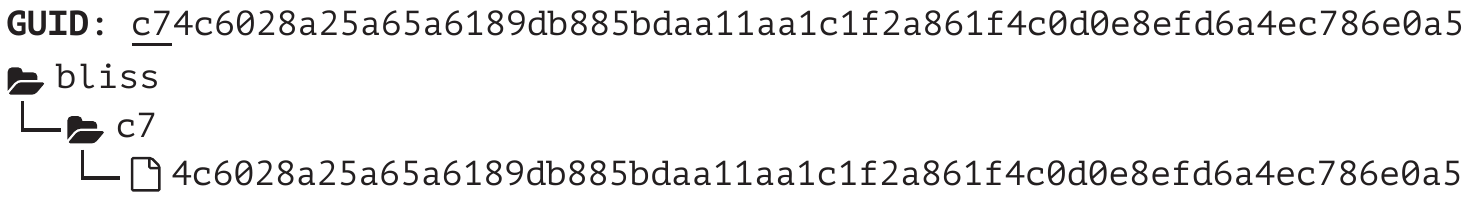}
\caption{Data repository: directory structure}
\end{figure}

Each data file contains a JSON object giving values of properties
corresponding to the classes it belongs to.
Note that any object in DiscreteZOO may belong to multiple classes --
there is no restriction that an object only belong to superclasses
of the most specific class.
For instance, a graph may simultaneously be
a cubic vertex-transitive graph and a split Praeger-Xu (SPX) graph,
with neither of these classes being a subclass of the other.
It is not required that every property be given --
those which have not yet been computed will be simply left out.

The data repository also contains information about datasets.
These are also stored as JSON files containing metadata,
a list of classes the dataset represents,
and a list of GUIDs of objects contained in the dataset.
Datasets may also be nested --
i.e., a parent dataset will contain all objects in the child dataset
(but will not necessarily describe the same object types).

The information stored in the data repository
can then be used to export the desired objects
to a form suitable for querying and accessing the objects,
e.g., an SQLite database.
Note that, due to performance reasons,
such a form may include additional identifiers
which are not guaranteed to be globally unique,
such as sequential row IDs.
This information is therefore not stored in the data repository,
which instead relies on GUIDs
for the purposes of referencing other objects in the database.

\subsubsection{Additional features of the SageMath package}

For the SageMath package,
it already suffices to describe a new collection
by adding the appropriate data and dataset files to the data repository.
Adding new object types, on the other hand, still require some programming.

The SQLite database in the SageMath package has a journaling system
that keeps track of the tables, rows and columns changed,
as well as of who introduced the changes.
This allows users to prepare contributions to the database.
To submit a contribution,
an author with an account on GitHub can make a pull request
to the data repository --
i.e., a request to the repository maintainers
to include the requested changes into the repository.
If accepted, the contribution is merged with the database,
the database downloads are updated
and other users can choose to update their local databases.

\subsection{Implementation and User Interfaces}

Each of the interfaces uses a simplified database optimised for its needs.
In addition to the core and package databases,
a user can download a subset of the core database for offline use with the SageMath package.
This local database can store any properties she computes.
From this, she can export the changes she makes to submit them to the core database.

\begin{figure}
\centering
\includegraphics[width=10cm,keepaspectratio]{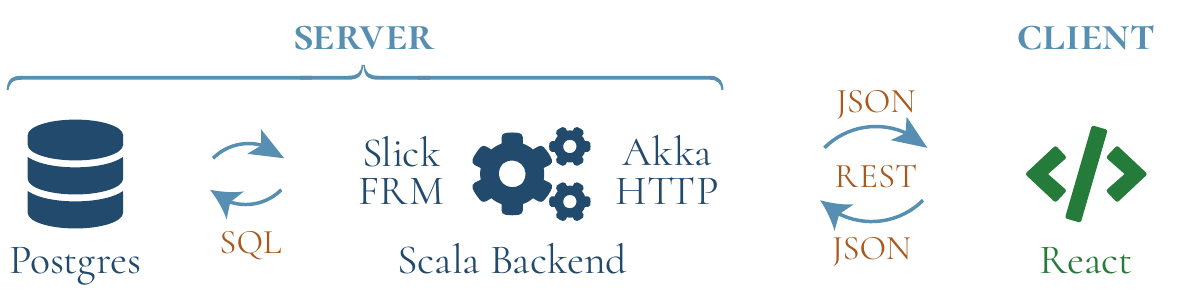}
\caption{Website stack}
\end{figure}

Both the website and the SageMath
package make it possible for the user
to search for objects and filter object sets.
The SageMath package supports
adding new objects and properties into the local database.  
If a researcher wants to submit some of this to the core database, 
the package helps with preparing the changes file.

\input{interfaces-web.tex}

\input{interfaces-sage.tex}

\section{Conclusions and Future Work}
\label{s:fw}

We built on Billey and Tenner's observations 
on the need of the mathematical community 
for searchable databases 
and described some further requirements
we consider important.
We observed that their work applies equally well
to databases of concrete mathematical objects,
which are similar to fingerprint databases of theorems,
but of slightly different flavour.
To get a sense of CMOs and databases 
we started an online catalogue, which
we hope to extend into a survey paper 
in the future.

DiscreteZOO is an experimental project 
towards a general purpose framework for CMOs.
The system is initially targeted at the community of researchers
working with symmetric graphs (and some related objects).
Above, we have already indicated
some short term goals and features for the future.
In particular, we would like to tackle coverage information~(\ref{re:coverage-information}).
An important practical consideration
will be to further automate updates to the database,
make them as declarative as possible
and thus decrease the effort needed to
contribute a new collection 
(\ref{re:collaborative}.\ref{re:declarative}).
First steps towards a ``programmable, mathematical API''
for mathematical knowledge bases
was investigated in~\cite{WieKohRab:vtuimkb17}
for the LMFDB.
We are going to explore the possibilities
offered by MMT~\cite{mmt:repo:on},
a framework for knowledge representation,
to support describing the CMOs and 
objects they contain.
In addition to the syntax, MMT also gives access
to a large body of formalised knowledge,
which would make implementing self-explanatory
aspects (\ref{re:self-explaining}) much easier.
We also plan to explore the
Math-in-the-Middle (MitM)~\cite{DehKohKon:iop16}
approach for interoperability~(\ref{re:interoperable}) with
computer algebra systems and
other databases.
MitM \emph{``allows translation between centrally formalized knowledge 
and the systems on the boundary via views and alignments''}.

\printbibliography

\end{document}

%% file: listings-python.tex
\definecolor{deepblue}{rgb}{0,0,0.5}
\definecolor{deepred}{rgb}{0.6,0,0}
\definecolor{deepgreen}{rgb}{0,0.5,0}
\definecolor{violet}{rgb}{0.5,0,0.5}

\DeclareFixedFont{\ttb}{T1}{txtt}{bx}{n}{8.5} 
\DeclareFixedFont{\ttm}{T1}{txtt}{m}{n}{8.5}  
\DeclareFixedFont{\ttc}{T1}{txtt}{m}{n}{8}  

\newcommand*{\mycommentstyle}[1]{%
  \begingroup
    \ttc
    \lstset{columns=fullflexible}%
    \color{deepgreen}%
    #1%
  \endgroup
}

\newcommand{\pythonstyle}{\lstset{
language=Python,
basicstyle=\small\ttm,
commentstyle=\mycommentstyle,
otherkeywords={self},             
keywordstyle=\ttb\color{deepblue},
emph={sage,__init__,True,False,None},  
emphstyle=\ttb\color{deepred},    
stringstyle=\color{violet},
frame=tb,                         
showstringspaces=false
}}

\lstnewenvironment{python}[1][]
{
\pythonstyle
\lstset{#1}
}
{}

\newcommand\pythoninline[1]{{\pythonstyle\lstinline!#1!}}

%% file: interfaces-web.tex
\subsection{Website\ednote{define API earlier, many screenshots}}
\label{ss:website}

The DiscreteZOO website is dedicated to simple searches,
downloads and displaying encyclopaedic information.

We will describe the website functionality on Use case~\ref{ex:partcubes}.
The census of cubic vertex-transitive graphs on up to $1280$ vertices~\cite{CVTCensus}
contains (among others) all graphs Tilen is interested in on up to $1280$ vertices.

\begin{figure}[h]
  \centering
  \includegraphics[width=\textwidth]{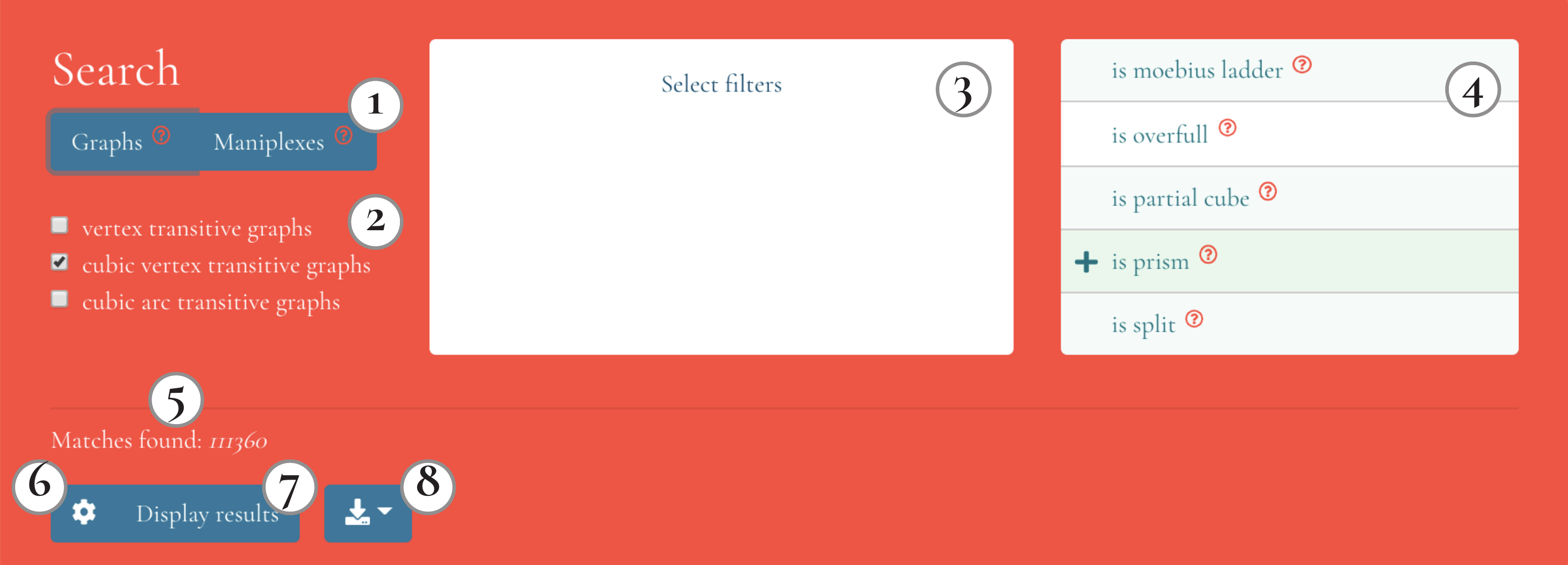}
  \caption{Search box}
  \label{figure:searchbox}
\end{figure}

\begin{wrapfigure}[7]{R}{0.45\textwidth}
\vspace*{-1.3em}
\centering
\includegraphics[width=0.45\textwidth,keepaspectratio]{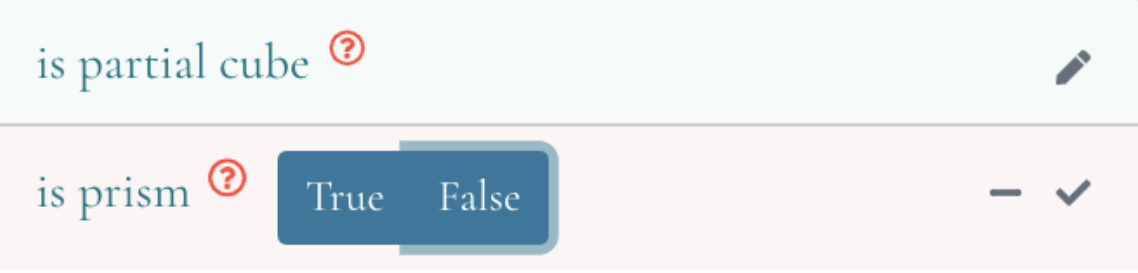}
  \caption{Selected filters}
  \label{figure:selectedfilters}
\end{wrapfigure}

\begin{enumerate}
\item In the search box (see Figure~\ref{figure:searchbox}),
he first chooses graphs as his objects~\circled{1}
\item and then only selects the collection of cubic vertex transitive graphs from the list
of graph collections~\circled{2}.
\item He first adds the filter \texttt{is partial cube}
from the available filters~\circled{4} to the selected filters~\circled{3},
chooses \texttt{True} and confirms the selection by clicking on the checkmark.
\item He knows that even prisms ($K_2 \square C_{2n}$)
are an infinite family of vertex-transitive cubic partial cubes 
and wants to remove them from the search results.
He first adds the filter \texttt{is prism}
from the available filters~\circled{4} to the selected filters~\circled{3},
chooses \texttt{False} to filter for those that are not prisms
and confirms selection by clicking on the checkmark
(see Figure~\ref{figure:selectedfilters}).
\item The number of matches~\circled{5} immediately updates with the number of objects in 
the database that match the new search criteria, every time
the chosen collections or filters get updated.
\item Finally, he downloads the list of graphs for SageMath
from the download dropdown menu~\circled{8}
and computes their isometric embeddings into hypercubes.
\end{enumerate}

To optimise responsiveness,
results are not displayed until the user presses
the ``Display results'' button (Figure~\ref{figure:searchbox}~\circled{6}).
Search result downloads do not require displaying the search results first.

The available filters area (Figure~\ref{figure:searchbox}~\circled{4}) contains all 
properties stored in the database for the chosen object type.
Numeric properties can be filtered to satisfy
simple equations (see Definition~\ref{def:numeric}).

\begin{figure}[h]
  \centering
  \includegraphics[width=\textwidth]{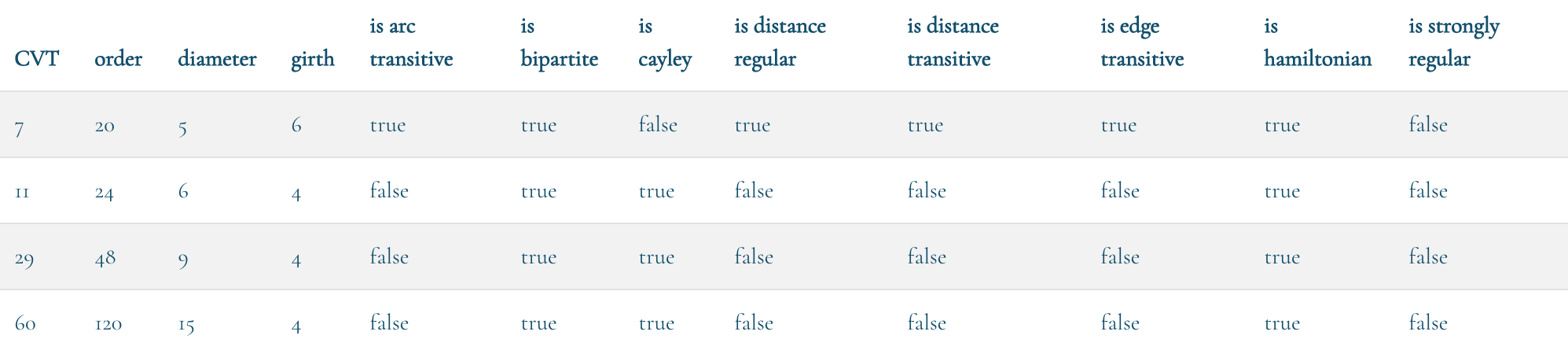}
  \caption{Displayed results}
  \label{figure:results}
\end{figure}

The following features are planned for implementation
in the near future.
\begin{enumerate}
\item Download
\begin{itemize}
\item graphs in the \texttt{sparse6}~\cite{GraphFormats} format as well as 
formats readable by computer algebra systems,
\item automorphism groups of symmetric objects in the search results
(in GAP and Magma code),
\item list of references in BibTeX (online resources, authors, papers)
relevant for the search results.
\end{itemize}
\item Copy data to clipboard (for example, 
graph in the \texttt{sparse6}, automorphism group code, references).
\item Tooltip definitions of concepts (object types, filters).
\end{enumerate}

%% file: interfaces-sage.tex
\subsection{SageMath package}
\label{ss:sage}

The objects in the database can also be accessed
using the SageMath package~\cite[{\tt DiscreteZOO-sage} repository]{ZooGroup}.
SageMath~\cite{Sage} is an open source computer algebra system
based on the Python programming language,
and it already provides many structures
for representing various mathematical objects.
The DiscreteZOO SageMath package defines its own structures
that inherit and override SageMath's structures,
thus allowing one to utilise the full potential of SageMath
while adding the functionality
of accessing precomputed properties in the database,
as well as storing newly computed properties back to the database
for later reuse.

After installing the package and the database,
one can either import the entire \pythoninline{discretezoo} package,
or load submodules as one needs them.
In the following examples,
we will use the submodule for cubic vertex-transitive graphs,
which includes the census of connected cubic vertex-transitive graphs
by Poto\v{c}nik, Spiga and Verret~\cite{CVTCensus}
(also known as the CVT census).
We only need to import the class \pythoninline{CVTGraph}
and the object \pythoninline{info}
from the \pythoninline{discretezoo.entities.cvt} submodule,
as well as the objects
from the \pythoninline{discretezoo.entities.cvt.fields} submodule.
The first submodule is intended for cubic vertex-transitive graphs,
while the latter contains the objects representing the precomputed properties
that one can use in search queries.


The \pythoninline{CVTGraph} class extends the \pythoninline{ZooGraph} class
representing general graphs in the database,
and the latter in turn extends both the general \pythoninline{ZooObject} class
and SageMath's \pythoninline{Graph} class.
One may construct a \pythoninline{CVTGraph} instance
by specifying the order and index as given
in the CVT census.
For instance, one may obtain the Petersen graph
and compare it to SageMath's builtin version
using SageMath's \pythoninline{is_isomorphic} method.

\begin{python}
sage: G = CVTGraph(10, 3)
sage: G.is_isomorphic(graphs.PetersenGraph())
True
\end{python}




A graph may also be constructed manually.
In this case, its canonical form is computed and a GUID is derived from it.
The GUID is then checked against the database,
and if there is a matching entry,
the precomputed properties will be loaded for the newly constructed graph.

\begin{python}
sage: CVTGraph([[(u, i) for u in GF(7) for i in (-1, 1)],
....:     lambda (u, i), (v, j): i != j and u*i + v*j in (1, 2, 4)])
Heawood graph: cubic vertex-transitive graph on 14 vertices, number 1
\end{python}

If the graph is not found in the database,
it will be checked for the required properties
(say, that it is cubic and vertex-transitive),
and added to the local copy of the database.
An entry indicating that a change to the database has occured
will then also be added,
thus allowing the user to track her changes
and later make a contribution to the main database.

In Use case~\ref{ex:partcubes},
Tilen is studying cubic vertex-transitive partial cubes.
He uses the \pythoninline{info} object to make queries to the database.
The \pythoninline{info.all} method returns a generator
yielding the requested graphs.
He restricts the queries by specifying conditions using the field objects.

\begin{python}
sage: gen = info.all(is_partial_cube, orderby=order) # order by the number of vertices
sage: next(gen) # first matching graph
3-Cube: cubic vertex-transitive graph on 8 vertices, number 2
sage: next(gen) # second matching graph
6-Prism: cubic vertex-transitive graph on 12 vertices, number 3
\end{python}

He notices that the small examples he has obtained are even prisms.
He now wants to obtain examples which are not prisms.

\begin{python}
sage: partial_cubes = list(info.all(is_partial_cube, ~is_prism))
sage: partial_cubes
[Desargues graph: cubic vertex-transitive graph on 20 vertices, number 7,
 Truncated Octahedron: cubic vertex-transitive graph on 24 vertices, number 11,
 Truncated Cuboctahedron: cubic vertex-transitive graph on 48 vertices, number 29,
 Truncated Icosidodecahedron: cubic vertex-transitive graph on 120 vertices, number 60]
\end{python}

He then examines the first graph in the list.
Asking whether the graph is a partial cube
will give the precomputed value which has been used to find it.

\begin{python}
sage: H = partial_cubes[0]
sage: H.is_partial_cube()
True
\end{python}

However, many methods in SageMath can also compute certificates
which can be used to verify the correctness of the result.
The certificate for a graph being a partial cube
is an isometric embedding into a hypercube --
i.e., a mapping from the vertices of the graph
to the vertices of the hypercube ($(0, 1)$-strings of a fixed length).
Since certificates are not stored in DiscreteZOO,
this triggers the computation.

\begin{python}
sage: _, cert = H.is_partial_cube(certificate=True)
sage: cert
{0: '00000', 1: '10111', 2: '11101', 3: '11011', 4: '11111', 5: '01011', 6: '11010', 7: '11100',
 8: '01101', 9: '10110', 10: '00111', 11: '11000', 12: '01001', 13: '10010', 14: '00011',
 15: '10100', 16: '00101', 17: '00010', 18: '00100', 19: '01000'}
\end{python}

Tilen thus sees that the Desargues graph
can be isometrically embedded into the $5$-dimensional hypercube.

On the other hand, in Use case~\ref{ex:smallgirth},
Primo\v{z} wants to count the number of graphs of each girth
in the CVT census.
He can do this with the \pythoninline{info.count} method.
Since the database contains some cubic vertex-transitive graphs
which only appear in other censuses,
he restricts his query to those graphs
for which the index from the census exists.

\begin{python}
sage: info.count(cvt_index) # number of graphs in the CVT census
111360
sage: info.count(cvt_index, groupby=girth) # break down by girth
{3: 160, 4: 5754, 5: 100, 6: 58674, 7: 192, 8: 13529, 9: 219, 10: 25806, 11: 80, 12: 5423, 13: 37,
 14: 1365, 15: 12, 16: 9}
\end{python}

He now wants to obtain a few small graphs of girth $7$ and diameter $4$.
This can be achieved using the \pythoninline{info.one} method,
which only returns a single graph.

\begin{python}
sage: info.one(girth == 7, diameter == 4, orderby=order)
Generalised Petersen graph (13, 5): cubic vertex-transitive graph on 26 vertices, number 5
sage: info.one(girth == 7, diameter == 4, orderby=order, offset=1) # the second such graph
Coxeter graph: cubic vertex-transitive graph on 28 vertices, number 6
\end{python}

The two main parts of the SageMath package implementation
are the \textbf{object classes} and the \textbf{database interfaces}.

The object class hierarchy closely follows the one
given in the class specifications repository (see Section~\ref{ss:repos}).
Upon initialisation, each class is initialised with the information
given in the corresponding JSON file.
Additionally, it also implements methods
needed for the functionality of the corresponding objects.
Thus, most generic functionality resides in the methods
of the \pythoninline{ZooEntity} and \pythoninline{ZooObject} classes.
The subclasses of \pythoninline{ZooObject}
also inherit their methods and properties from the appropriate classes
provided by SageMath.
These are however wrapped to deal with precomputed properties --
i.e., when a method is called without extra parameters,
a precomputed property will be returned if it is available,
otherwise it will be written back to the database once computed.
The methods such classes implement thus mostly deal with peculiarities
related to object construction
and calling certain methods inherited from SageMath.

To deal with queries,
each object class is accompanied by an \pythoninline{info} object
providing the user an interface for the queries,
and a \pythoninline{fields} submodule
containing objects representing the corresponding properties
(including the ones inherited from its superclasses).
These objects thus correspond to the properties of the query language
(see Section~\ref{sss:query}),
and they can be combined into more complex expressions
using arithmetic and logical operators.
Such expressions can then be given to the database interface,
which transforms it into a query for the appropriate database.
Currently, a generic SQL interface is available,
with SQLite and PostgreSQL database interfaces as its special cases.